\documentclass[11pt]{amsart}
\usepackage{float}
\usepackage[usenames]{color}
\usepackage{graphicx}
\usepackage{amssymb}
\usepackage{amscd}
\theoremstyle{plain}
\newtheorem{thm}{Theorem}[section]
\newtheorem{lemma}[thm]{Lemma}

\newtheorem{conj}[thm]{Conjecture}
\newtheorem{prop}[thm]{Proposition}

\theoremstyle{definition}

\newtheorem{example}[thm]{Example}

\newcommand{\fra}{\mathfrak{a}}
\newcommand{\frb}{\mathfrak{b}}

\newcommand{\frg}{\mathfrak{g}}
\newcommand{\frh}{\mathfrak{h}}
\newcommand{\fri}{\mathfrak{i}}
\newcommand{\frk}{\mathfrak{k}}
\newcommand{\frl}{\mathfrak{l}}

\newcommand{\frp}{\mathfrak{p}}

\newcommand{\frt}{\mathfrak{t}}
\newcommand{\fru}{\mathfrak{u}}

\newcommand{\bbC}{\mathbb{C}}

\newcommand{\bbN}{\mathbb{N}}

\newcommand{\be}{\begin {equation}}
\newcommand{\ee}{\end {equation}}
\newcommand{\bp}{\begin {proof}}
\newcommand{\ep}{\end {proof}}

\begin{document}

\title[Unitary Representations with Dirac for complex $E_6$]
{Unitary representations with non-zero\\ Dirac cohomology for complex $E_6$}

\author{Chao-Ping Dong}
\address[Dong]{Institute of Mathematics, Hunan University, Changsha 410082,
P.~R.~China}
\email{chaoping@hnu.edu.cn}
\thanks{Dong is supported by NSFC grant 11571097 and the China Scholarship Council.}

\abstract{This paper classifies the equivalence classes of irreducible unitary representations with nonvanishing Dirac cohomology for complex $E_6$. This is achieved by using our finiteness result, and by improving the computing method.}
\endabstract



\maketitle

\section{Introduction}
In his description of the wave function of the spin-$\frac{1}{2}$ particles, Dirac introduced the eponymous Dirac operator in \cite{Di} by using matrix algebra in 1928. This operator was a square root of the wave operator, and it led to the foundational Dirac equation in quantum mechanics. Mimicking the spirit of \cite{Di}, Parthasarathy introduced the geometric Dirac operator in representation theory of Lie groups \cite{P1} in 1972. This allowed him to construct most of the discrete series, and the construction was completed by Atiyah and Schmid \cite{AS}.

Let $G$ be a connected semisimple Lie group with finite center.
Let $\theta$ be the  Cartan involution of $G$ and assume that $K:=G^{\theta}$ is a maximal compact subgroup of $G$. Let $\frg=\frk\oplus\frp$ be the corresponding Cartan decomposition on the complexified Lie algebra level. Let  $U(\frg)$ be the universal enveloping algebra of $\frg$, and let $S$ be a spin module for the Clifford algebra $C(\frp)$. Let $\pi$ be any irreducible $(\frg, K)$ module. The Dirac operator $D$ lives in $U(\frg)\otimes C(\frp)$ and it  acts on $\pi\otimes S$. To understand the unitary dual of $G$ better, in 1997, Vogan formulated the notion of Dirac cohomology \cite{Vog97}, which  was defined to be $\widetilde{K}$-module
\begin{equation}\label{def-Dirac-cohomology}
H_D(\pi)=\text{Ker}\, D/ (\text{Im} \, D \cap \text{Ker} D).
\end{equation}
Here $\widetilde{K}$ is the spin double covering group of $K$. Vogan conjectured that whenever non-zero, Dirac cohomology  should  reveal  the infinitesimal character of the original module $\pi$.
This conjecture was proved by Huang and Pand\v zi\'c \cite{HP} in 2002 (see Theorem \ref{thm-HP}). Since then, Dirac cohomology became a new invariant for the study of Lie group representations.

We care most about the case that $\pi$ is unitary. Then $D$ is symmetric with respect to a natural inner product on $\pi\otimes S$, and
\begin{equation}\label{Dirac-unitary}
H_D(\pi)=\text{Ker}\, D=\text{Ker}\, D^2.
\end{equation}
 Parthasarathy's Dirac inequality \cite{P1,P2} now reads as that $D^2$ has non-negative eigenvalues on any $\widetilde{K}$-types of $\pi\otimes S$. Moreover, by Theorem \ref{thm-HP}, Dirac inequality becomes equality on some $\widetilde{K}$-types of $\pi\otimes S$ if and only if $H_D(\pi)$ is non-vanishing (see Proposition \ref{prop-D-spin-lowest} for more). Thus, among the entire unitary dual of $G$, those having non-zero Dirac cohomology are exactly the extreme ones.  Therefore,  classifying $\widehat{G}^{\mathrm{d}}$---the set of all the irreducible unitary representations (up to equivalence) with non-zero Dirac cohomology---should also be helpful for us to understand the entire unitary dual of $G$.

In this paper, we consider the special case that $G$ is a connected \emph{complex} simple Lie group. Recently, by analyzing Parthasarathy's Dirac inequality, by using results on cohomological induction mainly due to Vogan \cite{Vog84}, and by Theorem 6.1 of \cite{D}, we obtained in \cite{DD} a finiteness result:  $\widehat{G}^{\mathrm{d}}$ consists of finitely many scattered members (the scattered part) and finitely many strings of members (the string part), see Theorem \ref{thm-finite}. We also classified $\widehat{G}^{\mathrm{d}}$ for complex $F_4$ in \cite{DD}. Here, by improving the computing method, we report the following complete description of $\widehat{G}^{\mathrm{d}}$ for complex $E_6$.

\medskip
\noindent\textbf{Theorem A.}
\emph{The set $\widehat{E}_6^{\mathrm{d}}$ consists of $33$ scattered representations (see Table \ref{table-E6-scattered-part}) whose spin-lowest $K$-types are all unitarily small, and $213$ strings of representations (see Section \ref{sec-E6-string}).
Moreover, each representation $\pi\in\widehat{E}_6^{\mathrm{d}}$ has a unique spin-lowest $K$-type which occurs with multiplicity one.}
\medskip

In Theorem A, the notion unitarily small (\emph{u-small} for short) was introduced by Salamanca-Riba and Vogan in \cite{SV}, see Section \ref{sec-spin-norm}. The last statement of Theorem A  was motivated by \cite{H10}, where Huang kindly told the author that he  announced the following conjecture at a conference: each spin-lowest $K$-type of any $\pi\in\widehat{G}^{\mathrm{d}}$ should occur exactly once.

It is interesting to note that in the penultimate row of Table \ref{table-E6-scattered-part} sits the model representation due to McGovern \cite{Mc}, which is $K$-multiplicity free. A few other members there, say those described in Examples \ref{exam-comp}, \ref{exam-E6-2} and \ref{exam-E6-3}, may also be $K$-multiplicity free.

As deduced by Barbasch and Pand\v zi\'c on page 5 of \cite{BP} from Theorem \ref{thm-HP}, to find all the irreducible unitary representations with non-zero Dirac cohomology,  it suffices to consider the following candidates:
 \begin{equation}\label{BP}
J(\lambda, -s\lambda),
\end{equation}
where $s$ is an \emph{involution} in the Weyl group, and $\lambda$ is a weight such that $2\lambda$ is dominant integral and regular. Here $J(\lambda, -s\lambda)$ is the irreducible $(\frg, K)$ module with Zhelobenko parameters $\lambda_L=\lambda, \lambda_R=-s\lambda$, see Theorem \ref{thm-Zh}. At the end of Section 2.1, we will explain why the element $s$ in \eqref{BP} must be an involution.

 A little more thinking leads to the additional requirement that $\lambda-s\lambda$ should be a non-negative integer combination of simple roots, see \eqref{necc}. Surprisingly, all the calculations that we have carried out suggest that when put together, these  necessary conditions should become sufficient. Let us summarize this observation in the following.

\medskip
\noindent\textbf{Conjecture B.}
\emph{Let $G$ be a connected complex simple Lie group. The set $\widehat{G}^{\mathrm{d}}$ consists exactly of unitary representations $J(\lambda, -s\lambda)$, where $s$ is an involution, and $\lambda$ is a weight such that
\begin{itemize}
\item[$\bullet$] $2\lambda$ is dominant integral and regular;
\item[$\bullet$] $\lambda+s\lambda$ is an integral weight;
\item[$\bullet$] $\lambda-s\lambda$ is a non-negative integer combination of simple roots.
\end{itemize}
}
\medskip

The above conjecture holds for $A_1$-$A_6$, $B_2$-$B_4$, $C_2$-$C_4$, $D_4$-$D_6$, $G_2$, $F_4$ and $E_6$. Our calculations in type $A$ also lead to Conjecture \ref{conj-type-A}.

The paper is organized as follows.
We set up the notation and collect necessary preliminaries in Section 2. We discuss the automorphism $-w_0$ in Section 3, which will allow us to do reduction in calculations. Section 4 aims to improve the computing method of \cite{DD}. We figure out the scattered parts of $\widehat{G}^{\mathrm{d}}$ in Section 5 for some classical groups, and illustrate how to use this information to get the string part of $\widehat{E}_{6}^{\mathrm{d}}$ in Section 6. Finally, we determine the scattered part of $\widehat{E}_{6}^{\mathrm{d}}$ in Section 7.

\medskip
\emph{Acknowledgements.} I thank my advisor Prof.~Huang sincerely for sharing brilliant ideas with me during my PhD study. I also thank the math department of MIT for offering excellent working conditions. I am deeply grateful to the \texttt{atlas} mathematicians for many many things. Jian Ding had spent about one month to double-check all the calculations reported in this paper. In particular, the representation sitting in the $8$th row of Table \ref{table-E6-scattered-part} was originally missed. I thank him heartily for his time and carefulness.

\section{Preliminaries}

This section aims to set up the notation and collect some preliminaries.
Throughout this paper $\bbN=\{0, 1, 2, \dots\}$, $\mathbb{P}=\{1, 2, \dots\}$
and $\frac{1}{2}\mathbb{P}$ denotes the set of positive integers and positive half-integers.

Although some results in this section (say Theorem \ref{thm-HP} and Proposition \ref{prop-D-spin-lowest}) hold for real reductive Lie groups, for simplicity, we only quote them under the assumption that $G$ is a connected complex simple Lie group.
Let $\theta$ be the Cartan involution of $G$, and let $K:=G^{\theta}$ be a maximal compact subgroup of $G$.  Denote by
$\frg_0$ and $\frk_0$ the Lie algebras of $G$ and $K$, respectively. As usual, we drop the subscripts to denote the complexifications.
We denote by $\langle\, ,\,\rangle$  the Killing form form on $\frg$, which is negative definite on $\frk_0$ and positive definite on $\frp_0$. Moreover, $\frk$ and $\frp$ are orthogonal to each other under $\langle\, ,\,\rangle$. Let $\| \cdot\|$ be the norm corresponding to the Killing form.

Let $T$ be a maximal torus of $K$. Let $\fra_0=\sqrt{-1}\frt_0$
and $A=\exp(\fra_0)$. Then up to conjugation, $H=TA$ is the unique
$\theta$-stable Cartan subgroup of $G$.
We identify
\begin{equation}\label{identifications}
\frg\cong \frg_{0} \oplus
\frg_0, \quad
\frh\cong \frh_{0} \oplus
\frh_0, \quad \frt\cong \{(x,-x) : x\in
\frh_{0} \}, \quad \fra \cong\{(x, x) : x\in
\frh_{0} \}.
\end{equation}

Fix a Borel subgroup $B$ of $G$
containing $H$. Put $\Delta^{+}(\frg_0, \frh_0)=\Delta(\frb_0, \frh_0)$.
Then we have the corresponding simple roots $\alpha_1, \cdots, \alpha_l$ and fundamental weights $\varpi_1, \cdots, \varpi_l$. Set $[l]:=\{1,2, \dots, l\}$.  Denote by $s_i$ the simple reflection $s_{\alpha_i}$.
Let $\rho$ be the half sum of positive roots in $\Delta^{+}(\frg_0, \frh_0)$. In this paper, we always use the fundamental weights as a basis to express a weight. That is, $[n_1, \cdots, n_l]$ stands for the weight $\sum_{i=1}^{l} n_i \varpi_i$. For instance, $\rho=[1,1,1,1,1,1]$ for complex $E_6$. Set
$$
\Delta^+(\frg, \frh)=\Delta^+(\frg_0, \frh_0) \times \{0\} \cup \{0\} \times (-\Delta^+(\frg_0, \frh_0)).
$$
When restricted to $\frt$, we get $\Delta^+(\frg, \frt)$, $\Delta^+(\frk, \frt)$ and $\Delta^+(\frp, \frt)$. Denote by $\rho_c$ the half-sum of roots in $\Delta^+(\frk, \frt)$.
We denote by $W$ the Weyl group
$W(\frg_0, \frh_0)$, which has identity element $e$ and longest element $w_0$. Then $W(\frg, \frh)\simeq W \times W$.

\subsection{Zhelobenko classification}
The classification of irreducible admissible modules for complex Lie groups was obtained by Zhelobenko. Let $(\lambda_{L}, \lambda_{R})\in \frh_0^{*}\times
\frh_0^{*}$ be such that $\lambda_{L}-\lambda_{R}$ is
a weight of a finite dimensional holomorphic representation of $G$.
Using \eqref{identifications}, we can view $(\lambda_L, \lambda_R)$ as a real-linear functional on $\frh$, and write $\bbC_{(\lambda_L, \lambda_R)}$ as the character of $H$ with differential $(\lambda_L, \lambda_R)$ (which exists). Using \eqref{identifications} again, we have
$$
\bbC_{(\lambda_L, \lambda_R)}|_{T}=\bbC_{\lambda_L-\lambda_R}, \quad \bbC_{(\lambda_L, \lambda_R)}|_{A}=\bbC_{\lambda_L+\lambda_R}.
$$
Extend $\bbC_{(\lambda_L, \lambda_R)}$ to a character of $B$, and put $$X(\lambda_{L}, \lambda_{R})
:=K\mbox{-finite part of Ind}_{B}^{G}
(
\bbC_{(\lambda_L, \lambda_R)}
)
.$$

\begin{thm}\label{thm-Zh} {\rm (Zhelobenko \cite{Zh})}
The $K$-type with extremal weight $\lambda_{L}-\lambda_{R}$
occurs with multiplicity one in
$X(\lambda_{L}, \lambda_{R})$. Let
$J(\lambda_L,\lambda_R)$ be the unique subquotient of
$X(\lambda_{L}, \lambda_{R})$ containing this
$K$-type.
\begin{itemize}
\item[a)] Every irreducible admissible ($\frg$, $K$)-module is of the form $J(\lambda_L,\lambda_R)$.
\item[b)] Two such modules $J(\lambda_L,\lambda_R)$ and
$J(\lambda_L^{\prime},\lambda_R^{\prime})$ are equivalent if and
only if there exists $w\in W$ such that
$w\lambda_L=\lambda_L^{\prime}$ and $w\lambda_R=\lambda_R^{\prime}$.
\item[c)] $J(\lambda_L, \lambda_R)$ admits a nondegenerate Hermitian form if and only if there exists
$w\in W$ such that $w(\lambda_L-\lambda_R) =\lambda_L-\lambda_R , w(\lambda_L+\lambda_R) = -\overline{(\lambda_L+\lambda_R)}$.
\item[d)] The representation $X(\lambda_{L}, \lambda_{R})$ is tempered if and only if $\lambda_{L}+\lambda_{R}\in i\frh_0^*$. In this case,
$X(\lambda_{L}, \lambda_{R})=J(\lambda_{L}, \lambda_{R})$.
\end{itemize}
\end{thm}
Note that the $W\times W$ orbit of $(\lambda_L, \lambda_R)$ is the infinitesimal character of $J(\lambda_L, \lambda_R)$.
We call $\lambda_L$, $\lambda_R$ the\emph{ Zhelobenko parameters} for $J(\lambda_L, \lambda_R)$.
For instance, the trivial representation has $\lambda_L=\lambda_R=\rho$, while the model representation due to McGovern \cite{Mc} has $\lambda_L=\lambda_R=\rho/2$. We will also refer to $\lambda_L-\lambda_R$ (resp. $\lambda_L+\lambda_R$) as the \emph{$T$-parameter} (resp. \emph{$A$-parmeter}) of $J(\lambda_L, \lambda_R)$. The latter parameters are more convenient for the input of representations into \texttt{atlas}.

\begin{figure}[H]
\centering \scalebox{0.6}{\includegraphics{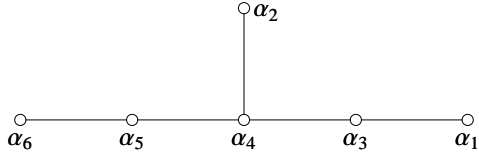}}
\caption{Dynkin diagram for $E_6$}
\end{figure}

\begin{example}\label{exam-atlas}
Let $G$ be complex $E_6$ (see Figure 1 for the labelling of the simple roots). Let
$$
s=s_4 s_5 s_6 s_5 s_1 s_3 s_2 s_4 s_1, \quad\lambda=[1, 1/2, 1/2, 1/2, 1/2, 1].
$$
Then $J(\lambda, -s\lambda)$ has $T$-parameter $[0, 4, 4, -4, 4, 0]$ and $A$-parameter $[2, -3, -3, 5, -3, 2]$. Performing the following commands allows us to input this representation into \texttt{atlas}:
\begin{verbatim}
set G=complex(simply_connected(E6))
set x=x(trivial(G))
set p=param(x, [0,4,4,-4,4,0,0,0,0,0,0,0],[2,-3,-3,5,-3,2,0,0,0,0,0,0])
\end{verbatim}
To test the unitarity of $J(\lambda, -s\lambda)$, we use the command
\begin{verbatim}
is_unitary(p)
\end{verbatim}
The output is
\begin{verbatim}
Value: true
\end{verbatim}
To look at the $K$-types of $J(\lambda, -s\lambda)$ up to the \texttt{atlas} height \texttt{h}, we use the command
\begin{verbatim}
branch_irr(p, h)
\end{verbatim}\hfill\qed
\end{example}

For convenience of reader, we repeat the explanation from \cite{BP} that the Weyl group element $s$ in \eqref{BP} must be an involution. Indeed, for $J(\lambda, -s\lambda)$ to be unitary, it should admit a non-degenerate Hermitian form in the first place. Thus by Theorem \ref{thm-Zh}(c), there exists $w\in W$ such that
$$
w(\lambda+s\lambda)=\lambda+s\lambda, \quad w(\lambda-s\lambda)=-\lambda+s\lambda.
$$
Therefore $w\lambda=s\lambda$ and $ws \lambda=\lambda$. Since $\lambda$ is regular, we have that $w=s$ and $ws=e$. Thus $s^2=e$, as desired.

\subsection{Dirac cohomology}
 Fix
an orthonormal basis $Z_1,\cdots, Z_n$ of $\frp_0$ with respect to
the inner product induced by $\langle\, ,\,\rangle$. Let $U(\frg)$ be the
universal enveloping algebra of $\frg$ and let $C(\frp)$ be the
Clifford algebra of $\frp$ with respect to $\langle\, ,\,\rangle$. The Dirac operator
$D\in U(\frg)\otimes C(\frp)$ is defined as
$$D=\sum_{i=1}^{n}\, Z_i \otimes Z_i.$$
It is easy to check that $D$ does not depend on the choice of the
orthonormal basis $Z_i$ and it is $K$-invariant for the diagonal
action of $K$ given by adjoint actions on both factors.

Let $\widetilde{K}$ be  the subgroup of $K\times \text{Spin}\,\frp_0$ consisting of all pairs $(k, s)$ such that $\text{Ad}(k)=p(s)$, where $\text{Ad}: K\rightarrow \text{SO}(\frp_0)$ is the adjoint action, and $p: \text{Spin}\,\frp_0\rightarrow \text{SO}(\frp_0)$ is the spin double covering map. Here $SO(\frp_0)$ is defined with respect to the Killing form on $\frp_0$.
If $\pi$ is a
($\frg$, $K$) module, and if $S$ denotes a spin module for
$C(\frp)$, then $\pi\otimes S$ is a $(U(\frg)\otimes C(\frp),
\widetilde{K})$ module. The action of $U(\frg)\otimes C(\frp)$ is
the obvious one, and $\widetilde{K}$ acts on both factors, on $\pi$
through $K$ and on $S$ through the spin group
$\text{Spin}\,{\frp_0}$.
Now the Dirac operator acts on $\pi\otimes S$, and the Dirac
cohomology of $\pi$ is the $\widetilde{K}$-module
defined in (\ref{def-Dirac-cohomology}).
By setting  the linear functionals on $\frt$ to be zero on $\fra$, we embed $\frt^{*}$ as a subspace of $\frh^{*}$. The following foundational result, conjectured
by Vogan,  was proved by Huang and Pand\v zi\'c \cite{HP}.

\begin{thm}{\rm (Huang and Pand\v zi\'c)}\label{thm-HP}
Let $\pi$ be an irreducible ($\frg$, $K$) module.
Assume that the Dirac
cohomology of $\pi$ is nonzero, and that it contains the $\widetilde{K}$-type $E_{\gamma}$ with highest weight $\gamma\in\frt^{*}\subset\frh^{*}$. Then the infinitesimal character of $\pi$ is conjugate to
$\gamma+\rho_{c}$ under $W(\frg,\frh)$.
\end{thm}

\subsection{Spin norm and spin lowest $K$-type}\label{sec-spin-norm}
The notions of spin norm and spin-lowest $K$-type were introduced in the author's thesis for real reductive Lie groups. They are motivated for the classification of irreducible unitary representations with non-zero Dirac cohomology. Let us recall them for complex Lie groups.
We identify a $K$-type $\delta$ with its highest weight. Then
\begin{equation}\label{Spin-norm-K-type}
\|\delta\|_{\mathrm{spin}}= \| \{\delta-\rho\} + \rho \|
\end{equation}
is the \emph{spin norm} of the $K$-type $\delta$.
Here $\{\delta-\rho\}$ is the unique dominant weight to which $\delta-\rho$ is conjugated under the action of $W$. Recall that $\delta$ is u-small in the sense of Salamanca-Riba and Vogan \cite{SV} if and only if $\delta$ lies in the convex hull of the $W$-orbit of $2\rho$. In such a case, by Lemma 2.3 of \cite{D2},
$$
\|\rho\|\leq \|\delta\|_{\mathrm{spin}}\leq \|2\rho\|.
$$
For any
irreducible admissible ($\frg$, $K$)-module $\pi$, we define
\begin{equation}\label{spin-norm-X}
\|\pi\|_{\mathrm{spin}}=\min \|\delta\|_{\mathrm{spin}},
\end{equation}
where $\delta$ runs over all the $K$-types occurring in $\pi$.
We call $\delta$ a \emph{spin lowest $K$-type} of $\pi$ if
it occurs in $\pi$ and
$\|\delta\|_{\mathrm{spin}}=\|\pi\|_{\mathrm{spin}}$.

Let us recall Proposition 3.3 of \cite{D} for complex Lie groups. It is a combination of the ideas and results of
Parthasarathy \cite{P1, P2}, Vogan \cite{Vog97}, Huang and Pand\v
zi\'c \cite{HP}.

\begin{prop}\label{prop-D-spin-lowest}
For any irreducible unitary ($\frg$, $K$)-module $\pi$ with
infinitesimal character $\Lambda$, let $\delta$ be any $K$-type
occurring in $\pi$. Then
\begin{enumerate}
\item[a)] $\|\pi\|_{\mathrm{spin}}\geq\|\Lambda\|$, and the equality holds if and only if $H_D(\pi)$ is non-zero.
\item[b)] $\|\delta\|_{\mathrm{spin}}\geq \|\Lambda\|$, and the equality holds
if and only if $\delta$ contributes to $H_D(\pi)$.
\item[c)] If $H_D(\pi)\neq 0$, it is exactly the spin lowest $K$-types of $\pi$
that contribute to $H_D(\pi)$.
\end{enumerate}
\end{prop}
In view of the above proposition, spin norm and spin lowest $K$-type give the right framework for the classification of $\widehat{G}^{\mathrm{d}}$.

\subsection{Vogan pencil}
Let $\beta$ be the highest root.
The following result is a special case of Lemma 3.4 and Corollary 3.5 of \cite{Vog81}. It coarsely describes the $K$-type pattern for an infinite-dimensional irreducible ($\frg$, $K$)-module $\pi$.

\begin{prop}\label{Vogan-K-types-pattern}{\rm (Vogan)}
Let $G$ be a connected complex simple Lie group. Then for any infinite-dimensional irreducible ($\frg$, $K$)-module $\pi$, there is a unique set
$$\{\mu_{\fri} \,|\, \fri \in  I \} \subseteq
i \frt_{0}^{*}$$ of  dominant integral weights such that all the
$K$-types of $\pi$ are precisely
$$\{\mu_{\fri} +n \beta
\,|\, \fri \in  I, n\in\bbN \}.$$
\end{prop}

We call a set of $K$-types
\begin{equation}\label{P-delta}
P(\delta):=\{\delta +n \beta \,|\, n\in\bbN \}
\end{equation}
a \emph{Vogan pencil}. For instance, $P(0)$ denotes the pencil starting from the trivial $K$-type. We also set
 \begin{equation}\label{P-mu-prime}
P_{\delta}:=\min \{\|\delta +n \beta \|_{\mathrm{spin}}\,|\, n\in\bbN \}.
\end{equation}
Calculating $P_\delta$ will be a vital step in our computing method in Section \ref{sec-comp}. On this aspect, we mention that by Theorem 1.1 of \cite{D2},
\begin{equation}\label{P-mu}
P_\delta=
\begin{cases}
\min \{\|\delta +n \beta \|_{\mathrm{spin}}\,|\, \delta+n\beta \mbox{ is u-small}\} & \mbox{ if $\delta$ is u-small};\\
\|\delta\|_{\mathrm{spin}} & \mbox{ otherwise.}
\end{cases}
\end{equation}

\subsection{A necessary condition}
As mentioned in the introduction, to find all the irreducible unitary representations with non-zero Dirac cohomology,  it suffices to consider the candidates in \eqref{BP}. We can add one more requirement here. Indeed, suppose that $J(\lambda, -s\lambda)$ is a member of  $\widehat{G}^{\mathrm{d}}$. Then by Theorem \ref{thm-HP} and Proposition \ref{prop-D-spin-lowest}, it has a spin lowest $K$-type $\delta$ such that
$$
\{\delta-\rho\}+\rho=2\lambda.
$$
Since $\{\delta-\rho\}=\delta-\rho+\sum_i p_i \alpha_i$ for some $p_i\in \bbN$, it follows that
\begin{equation}\label{necc-1}
2\lambda=\delta+\sum_i p_i \alpha_i.
\end{equation}
On the other hand, put $\mu:=\{\lambda+s\lambda\}$. Then
\begin{equation}\label{necc-2}
\mu=\lambda+s\lambda+\sum_i q_i \alpha_i,
\end{equation}
where $q_i\in\bbN$. Since $\mu$ is the lowest $K$-type of $J(\lambda, -s\lambda)$, by Frobenius reciprocity and the highest weight theorem, we have
\begin{equation}\label{necc-3}
\delta=\mu+\sum_i r_i\alpha_i,
\end{equation}
for some $r_i\in\bbN$. Combining \eqref{necc-1}, \eqref{necc-2}, and \eqref{necc-3} gives
\begin{equation}\label{necc}
\lambda-s\lambda=\sum_i n_i \alpha_i,
\end{equation}
where $n_i=p_i+q_i+r_i\in\bbN$.

\subsection{A description of $\widehat{G}^{\mathrm{d}}$}
 One key idea in \cite{DD} was to arrange the representations \eqref{BP} into $s$-families.  More precisely, let us denote
\begin{equation}\label{Lambda-s}
 \Lambda(s):=\left\{\lambda=[\lambda_1, \dots, \lambda_l]\mid 2\lambda_i\in\mathbb{P}, \lambda+s\lambda \mbox{ is integral, and } \lambda-s\lambda \mbox{ satisfies } \eqref{necc}\right\}.
\end{equation}
We call $\Lambda(s)$ and the corresponding representations $J(\lambda, -s\lambda)$ an \emph{$s$-family}. Note that an $s$-family has infinitely many members. For instance, the $e$-family consists of tempered representations, and they are handled in Section 4 of \cite{DD}; while on the other extreme, spherical representations live in the $w_0$-family, and they are considered in Section 5 of \cite{DD}.

Let $I$ be a \emph{non-empty} subset of $[l]$.  We call
\begin{equation}\label{I-string}
 \{\lambda\in\Lambda(s)\mid \lambda_i \mbox{ varies for } i\in I \mbox{ and } \lambda_j \mbox{ is fixed for } j\in[l]\setminus I\}
\end{equation}
and the corresponding representations $J(\lambda, -s\lambda)$ an \emph{$(s, I)$-string}.
When $s$ is known from the context, we may call it an \emph{$I$-string} or just a \emph{string}.

Fix an involution $s$. Put
\begin{equation}\label{I-s}
 I(s):=\left\{i\in[l]\mid s(\varpi_i)=\varpi_i\right\}.
\end{equation}
As shown in Lemma 3.1 of \cite{DD}, the set $I(s)$ consists of the indices $i$ such that the simple reflection $s_i$ does not appear in one reduced expression of $s$.

\begin{thm}\label{thm-finite} \emph{(Theorem A of \cite{DD})\footnote {This result has been partially generalized to real reductive Lie groups in \cite{D17}.}}
Fix an involution $s$ of a complex connected simple Lie group $G$.
If $I(s)$ is empty, then the $s$-family contains at most finitely many members of $\widehat{G}^{\mathrm{d}}$. If $I(s)$ is non-empty, then the $s$-family contains at most finitely many $I(s)$-strings of members of $\widehat{G}^{\mathrm{d}}$. The latter representations are cohomologically induced from members of $\widehat{L}_s^{\mathrm{d}}$ sitting in the $s$-family of $L_s$, and they are all in the good range. Here $L_s\supseteq HA$ is the $\theta$-stable Levi subgroup of $G$ corresponding to the simple roots $\{\alpha_i\mid i\notin I(s)\}$.
\end{thm}

We call the members of $\widehat{G}^{\mathrm{d}}$ coming from those $s$-families such that $I(s)$ are empty the \emph{scattered part} of $\widehat{G}^{\mathrm{d}}$, while  we call the remaining members  of $\widehat{G}^{\mathrm{d}}$ the \emph{string part}. Recall that by the proof of Proposition 3.5 of \cite{DD}, $J(\lambda, -s\lambda)$ is cohomologically induced from the irreducible $(\frl_s, L_s \cap K)$ module with Zhelobenko parameters $\lambda -\rho(\fru_s)/2$ and $-s(\lambda-\rho(\fru_s)/2)$. Here $P_s$ is the $\theta$-stable parabolic subgroup of $G$ with Levi factor $L_s$, $\frp_s=\frl_s +\fru_s$ is the Levi decomposition of the complexified Lie algebbra of $P_s$, and $\rho(\fru_s)$ is the half-sum of the positive roots in $\Delta(\fru_s, \frh)$. The induction is always in the good range. Thus it preserves unitarity \cite{Vog84} and the non-vanishing of Dirac cohomology \cite{D}.  Therefore, to figure out $\widehat{G}^{\mathrm{d}}$, it suffices to pin down the scattered parts of $\widehat{G}^{\mathrm{d}}$ and  $\widehat{L}_{\mathrm{ss}}^{\mathrm{d}}$. Here $L_{\mathrm{ss}}=[L_s, L_s]$.
Moreover, by Proposition 3.4 of \cite{DD}, it suffices to consider finitely many candidates representations to determine the scattered part of $\widehat{G}^{\mathrm{d}}$. An explicit algorithm will be presented in Section \ref{sec-comp} to sieve out these finite candidates. To sum up, after a finite calculation, one can completely determine $\widehat{G}^d$.

\section{The automorphism $-w_0$}\label{sec-w0}
In this section, we assume that $-w_0\neq 1$, which happens exactly when $G$ is $A_n$ ($n\geq 2$), $D_{2n+1}$ and $E_6$. As we shall see, the map $-w_0$ gives an automorphism of the complex Lie group $G$ and will allow us to do reduction in studying representations of $G$.

Let $s\in W$ be an involution. Then it is obvious that $w_0 s w_0$ is still an involution. If $w_0 s w_0=s$, we say the involution $s$ is \emph{self-dual}; otherwise, $s^\prime:=w_0 s w_0$ is another involution. In the latter case, we say that $s$ and $s^\prime$ are \emph{dual} to each other. For instance, $E_6$ has $892$ involutions in total, among which $140$ are self-dual.

Since $w_0\rho=-\rho$, the following lemma is immediate.
\begin{lemma}\label{lemma-involution-dual}
 We have the following.
\begin{itemize}
\item[a)] Two involutions $s$ and $s^\prime$ are dual to each other if and only if
$(-w_0) s\rho=s^\prime\rho$.
\item[b)] The involution $s$ is self-dual if and only if
$(-w_0) s\rho=s\rho$.
\end{itemize}
\end{lemma}

Under our assumption that $-w_0\neq 1$, we have that $-w_0$ is not an element of $W$. Therefore, in view of Theorem \ref{thm-Zh}, the two representations $J(\lambda_L, \lambda_R)$ and $J(-w_0\lambda_L, -w_0\lambda_R)$ are inequivalent. However,
  since $-w_0$ gives an automorphism of $G$, they share the same unitarity, while the dual $K$-type pattern. Thus we say they are \emph{dual} to each other, and use $\sim$ to denote this relation. That is,
  \begin{equation}
  J(\lambda_L, \lambda_R)\sim J(-w_0\lambda_L, -w_0\lambda_R).
  \end{equation}
In such a case, we can \emph{fold} the two representations by studying only one of them. This will reduce the work load.

Now suppose that $s$ and $s^\prime$ are dual to each other. Then we have
$$
J(\lambda, -s \lambda)\sim J(-w_0\lambda, w_0s \lambda)=J(-w_0\lambda, (-s^\prime)(-w_0)\lambda).
$$
Thus it suffices to study the $s$-family. On the other hand, if $s$ is self-dual, we have
$$
J(\lambda, -s \lambda)\sim J(-w_0\lambda, w_0 s \lambda)=J(-w_0\lambda, (-s)(-w_0)\lambda).
$$
Therefore, within the $s$-family, whenever $\lambda\neq -w_0\lambda$, it suffices to consider the parameter $\lambda$.

In the following sections, we always present the \emph{folded version} of the scattered part of $\widehat{G}^{\mathrm{d}}$: We mark an involution $s$ with a star whenever it is not self-dual. When $s$ is self-dual, we mark the parameter $\lambda$ with a star whenever $\lambda\neq -w_0\lambda$. In other words, the appearance of a star always indicates the existence of two representations which are dual to each other. Thus we can unfold and restore the entire scattered part easily.

\section{The improved computing method}\label{sec-comp}
This section aims to introduce a method that allows us to compute all the members of $\widehat{G}^{\mathrm{d}}$ in any $s$-family such that $I(s)$ is empty.  We proceed as follows:
\begin{itemize}
\item[$\bullet$] collect the finitely many $\lambda\in \Lambda(s)$ such that $\lambda-s\lambda=\sum_i n_i\alpha_i$, where $n_i\in\bbN$, and that
\begin{equation}\label{step-1}
\|\lambda-s\lambda\|^2\leq \|2\rho\|^2.
\end{equation}
\item[$\bullet$] further collect from the previous step those $\lambda$ satisfying
\begin{equation}\label{step-2}
\|2\lambda\|^2 \leq P_{\mu}^2,
\end{equation}
where $\mu:=\{\lambda+s\lambda\}$ is the lowest $K$-type of $J(\lambda, -s\lambda)$ and $P_{\mu}$ is defined in \eqref{P-mu}.
\item[$\bullet$] For the remaining $\lambda$, use  \texttt{atlas} \cite{ALTV,At} to study the unitarity and $K$-types of $J(\lambda, -s\lambda)$.
\end{itemize}

Let us explain why the method works. For the first step, as deduced in Section 3 of \cite{DD},
$$
\|2\lambda\|^2-\|\mu\|_{\mathrm{spin}}^2=\|\lambda-s\lambda\|^2-g(\lambda),
$$
where $g(\lambda):=2\langle \{\mu-\rho\}-(\mu-\rho), \rho\rangle\leq\|2\rho\|^2$, see Lemma 3.3 of \cite{DD}. Thus, if $\lambda$ does not meet the requirement \eqref{step-1}, we would have that
$$\Delta_1(\lambda):=\|2\lambda\|^2-\|\mu\|_{\mathrm{spin}}^2>0.
$$ Therefore, the corresponding representation $J(\lambda, -s\lambda)$ is non-unitary by Dirac inequality. Since $I(s)$ is assumed to be empty, by Lemma 3.2 of \cite{DD}, $\|\lambda-s\lambda\|^2$ is a homogeneous quadratic polynomial in terms of $\lambda_i$. Moreover, each term $\lambda_i^2$ has a positive coefficient, while each term $\lambda_i \lambda_j$ ($i<j$) has a nonnegative coefficient. Thus there are finitely many $\lambda$ satisfying \eqref{step-1}. For the second step, if $\lambda$ does not meet the condition \eqref{step-2}, we would have that
$$\Delta_2(\lambda):=\|2\lambda\|^2-P_{\mu}^2>0.
$$
Thus, again by Dirac inequality, the corresponding representation $J(\lambda, -s\lambda)$ is non-unitary.

The current method improves the previous one \cite{DD} mainly at the first step. Indeed, most of our energy in sieving out the candidate representations for $\widehat{F}_4^{\mathrm{d}}$ was spent in obtaining the specific values of $g(\lambda)$ via case-by-case analysis, see Section 8 of \cite{DD}.
To carry out a similar analysis for type $E$ is next to impossible. This motivates our first step: by adopting the uniform bound $\|2\rho\|^2$ of the function $g(\lambda)$, we no longer need to do case-by-case analysis of its values. For the second step, we note that using the distribution of the spin norm along the Vogan pencil $P(\mu)$ is very efficient in practical calculation, see Example \ref{exam-comp}.

Recall that in Section 5 of \cite{DD}, which is essentially known in Section 7 of \cite{D}, we have an effective  way to deal with the spherical unitary dual living in the $w_0$-family.  The current method extends it to any $s$-family such that $I(s)$ is non-empty. The extension is non-trivial in the following sense: unlike the spherical representations, now the lowest $K$-type $\mu:=\{\lambda+s\lambda\}$ of $J(\lambda, -s\lambda)$ varies according to $\lambda$ and $s$. The key ingredient leading to the extension is the analysis of Parthasarathy's Dirac inequality carried out in Section 3 of \cite{DD}.

Remark that we used \texttt{Mathematica} to carry out the first two steps, and the pdf version was uploaded on \texttt{ReserachGate} via the link
\begin{verbatim}
https://www.researchgate.net/publication/320110729_E6-s-family-genuine
\end{verbatim}
We gave explanations to the codes.  Thus the reader can pick up them easily. An interested reader may also modify these codes to investigate other complex Lie groups. On the other hand, one can carry out the third step using the \texttt{atlas} commands in Example \ref{exam-atlas}.

\begin{example}\label{exam-comp}
Let $G$ be complex $E_6$. Consider the self-dual involution
$$
s=s_4 s_5 s_6 s_5 s_1 s_3 s_2 s_4 s_1.
$$
Note that $s\rho=[-2,5,6,-7,6,-2]$ (recall Lemma \ref{lemma-involution-dual}).

The first step leaves us with $124048$ candidate representations. However, after carrying out the second step, only the following three $\lambda$ survive:
$$
[1/2, 1/2, 1, 1/2, 1/2, 1], [1, 1/2, 1/2, 1/2, 1/2, 1], [1, 1/2, 1/2,
1/2, 1, 1/2].
$$
The first and the third are dual to each other. Thus we can fold them by omitting the third one:
$$
[1/2, 1/2, 1, 1/2, 1/2, 1], [1, 1/2, 1/2, 1/2, 1/2, 1].
$$
Then by \texttt{atlas}, we know that only the second $\lambda$ gives a unitary representation $J(\lambda, -s\lambda)$, which has $T$-parameter $[0, 4, 4, -4, 4, 0]$ and $A$-parameter $[2, -3, -3, 5, -3, 2]$. (Recall Example \ref{exam-atlas}.) This representation has a unique spin-lowest $K$-type $[1,1,0,3,0,1]$ which occurs once.   Moreover,
 $$
 \|[1,1,0,3,0,1]\|_{\mathrm{spin}}=\|2\lambda\|.
 $$
Thus it is a member of $\widehat{E}_6^{\mathrm{d}}$ by Proposition \ref{prop-D-spin-lowest}.
This representation sits in the first row of Table \ref{table-E6-scattered-part}. \hfill\qed
\end{example}

\section{The scattered part of $\widehat{G}^{\mathrm{d}}$ for some classical groups}
\label{sec-classical}
This section aims to describe the scattered part of $\widehat{G}^{\mathrm{d}}$ for some classical groups with small ranks. This information will be needed later to form the string part of $\widehat{E}_6^{\mathrm{d}}$. For convenience, in each table, we always present the \emph{folded version}. (Recall the last paragraph of Section 3.) Therefore, one just needs to pay attention to each star to restore the entire scattered part. In particular, in each table, $N_G$ equals the number of rows plus the number of stars. Here $N_G$ denotes the cardinality of the scattered part of $\widehat{G}^{\mathrm{d}}$. We note also that the trivial representation always sits in the last row of each table.

\subsection{The scattered part of $\widehat{A}_{i}^{\mathrm{d}}$ ($1\leq i\leq 5$)}
One can calculate that the scattered part of $\widehat{A}_1^{\mathrm{d}}$ consists of the trivial representation, and that of $\widehat{A}_2^{\mathrm{d}}$ consists of the trivial representation and the model representation. We list the folded version of scattered parts of $\widehat{A}_{i}^{\mathrm{d}}$ ($3\leq i\leq 5$) in Tables \ref{table-A3-scattered-part}--\ref{table-A5-scattered-part}. To sum up, we have $N_{A_1}=1$, $N_{A_2}=2$, $N_{A_3}=4$, $N_{A_4}=8$ and $N_{A_5}=16$.

One can also calculate that $N_{A_6}=32$. Thus we make the following.
\begin{conj}\label{conj-type-A}
We have $N_{A_n}=2^{n-1}$.
\end{conj}

\begin{table}
\centering
\caption{The scattered part of $\widehat{A}_3^{\mathrm{d}}$ (folded version)}
\begin{tabular}{l|c|c|c|r}
$s\rho$ &   $\lambda$   & spin LKT & mult &  u-small\\
\hline
$[1,-3,1]$ & $\rho/2$ & $\rho$  & $1$ & Yes\\
$[-2,1,-2]$ & $[1/2, 1/2, 1]^*$ & $[2, 0, 1]$  & $1$ & Yes\\
$-\rho$ & $\rho$ & $[0,0,0]$ & $1$ & Yes
\end{tabular}
\label{table-A3-scattered-part}
\end{table}

\begin{table}
\centering
\caption{The scattered part of $\widehat{A}_4^{\mathrm{d}}$ (folded version)}
\begin{tabular}{l|c|c|c|r}
$s\rho$ &   $\lambda$   & spin LKT & mult &  u-small\\
\hline
$[-2,3,-4,2]^*$ & $[1, 1/2, 1/2, 1/2]$ & $[1, 0, 2, 1]$  & $1$ & Yes\\
$[-3,1,1,-3]$ & $\rho/2$ & $\rho$  & $1$ & Yes\\
$[-2,-1,2,-3]^*$ & $[1,1,1/2,1/2]$ & $[1,0,0,3]$   & $1$ & Yes\\
$-\rho$ & $\rho/2$ & $\rho$  & $1$ & Yes\\
$-\rho$ & $[1,1/2,1/2,1]$ & $2\beta$ & $1$ & Yes\\
$-\rho$ & $\rho$ & $[0,0,0,0]$  & $1$ & Yes
\end{tabular}
\label{table-A4-scattered-part}
\end{table}

\begin{table}
\centering
\caption{The scattered part of $\widehat{A}_5^{\mathrm{d}}$ (folded version)}
\begin{tabular}{l|c|c|c|r}
$s\rho$ &   $\lambda$   & spin LKT & mult &  u-small\\
\hline
$[-3,1,3,-5,3]^*$ & $\rho/2$ & $\rho$  & $1$ & Yes\\
$[-2,4,-5,4,-2]$ & $[1, 1/2, 1/2, 1/2, 1]$ & $[1, 0, 3, 0, 1]$  & $1$ & Yes\\
$[-2,-1,4,-5,3]^*$ & $[1, 1, 1/2, 1/2, 1/2]$ & $[1, 0, 0, 3, 1]$  & $1$ & Yes\\
$[-4,2,-1,2,-4]$ & $[1/2, 1/2, 1, 1/2, 1/2]$ & $[2, 1, 0, 1, 2]$  & $1$ & Yes\\
$[-2,-1,-1,3,-4]^*$ & $[1, 1/2, 1/2, 1/2, 1/2]$ & $[1, 1, 0, 2, 1]$  & $1$ & Yes\\
$[-2,-1,-1,3,-4]^*$ & $[1, 1, 1, 1/2, 1/2]$ & $[1, 0, 0, 0, 4]$  & $1$ & Yes\\
$[1,-3,1,-3,1]$ & $\rho/2$ & $\rho$  & $1$ & Yes\\
$[-2,1,-3,1,-2]$ & $[1/2, 1/2, 1/2, 1/2, 1]^*$ & $[1, 2, 0, 1, 1]$  & $1$ & Yes\\
$[-1,-2,1,-2,-1]$ & $[1, 1/2, 1/2, 1, 1]^*$ & $[3, 0, 0, 0, 2]$  & $1$ & Yes\\
$-\rho$ & $\rho$ & $[0,0,0,0,0]$  & $1$ & Yes
\end{tabular}
\label{table-A5-scattered-part}
\end{table}

\subsection{The scattered part of $\widehat{D}_{4}^{\mathrm{d}}$}

Note that $-w_0=1$ for $D_4$. However, we can use the automorphism which interchanges $\alpha_3$ and $\alpha_4$ while preserving $\alpha_1$ and $\alpha_2$ to play the role of $-w_0$. In this sense, we present the folded version of the scattered part of $\widehat{D}_{4}^{\mathrm{d}}$ in Table \ref{table-D4-scattered-part}. To keep the folded versions in a unified style, we do not use other automorphisms of $D_4$. Note that the penultimate row is the model representation due to McGovern \cite{Mc}, and that $N_{D_4}=9$.

\begin{table}
\centering
\caption{The scattered part of $\widehat{D}_4^{\mathrm{d}}$ (folded version)}
\begin{tabular}{l|c|c|c|r}
$s\rho$ &   $\lambda$   & spin LKT & mult &  u-small\\
\hline
$[3,-5,3,3]$ & $\rho/2$ & $\rho$  & $1$ & Yes\\
$[-5,3,-1,-1]$ & $[1/2, 1/2, 1, 1]$ & $[3,1,0,0]$  & $1$  & Yes \\
$[-1,3,-5,-1]^*$ & $[1, 1/2, 1/2, 1]$ & $[0,1,3,0]$  & $1$ & Yes \\
$[-1,-3,1,1]$ & $[1, 1/2, 1/2, 1/2]$ & $[1,2,0,0]$  & $1$  & Yes \\
$[1,-3,-1,1]^*$ & $[1/2, 1/2, 1, 1/2]$ &  $[2,1,0,0]$ & $1$  & Yes \\
$-\rho$ & $\rho/2$ & $\rho$  & $1$ & Yes\\
$-\rho$ & $\rho$ & $[0,0,0,0]$  & $1$ & Yes
\end{tabular}
\label{table-D4-scattered-part}
\end{table}

\subsection{The scattered part of $\widehat{D}_{5}^{\mathrm{d}}$}
In this case, $-w_0$ interchanges $\alpha_4$ and $\alpha_5$, while preserving other simple roots.
The information is presented in Table \ref{table-D5-scattered-part}. We remark that the penultimate row is the model representation, and that $N_{D_5}=17$.

\begin{table}
\centering
\caption{The scattered part of $\widehat{D}_5^{\mathrm{d}}$ (folded version)}
\begin{tabular}{l|c|c|c|r}
$s\rho$ &   $\lambda$   & spin LKT & mult &  u-small\\
\hline
$[-2,5,-6,4,4]$ & $[1, 1/2, 1/2, 1/2, 1/2]$ & $[1, 0, 2, 1, 1]$  & $1$ & Yes\\
$[5,-7,5,-1,-1]$ & $[1/2, 1/2, 1/2, 1, 1]$ & $[1, 3, 1, 0, 0]$  & $1$ & Yes\\
$[-4,3,4,-6,-2]^*$ & $[1/2, 1/2, 1/2, 1/2, 1]$ & $[2, 0, 1, 2, 1]$  & $1$ & Yes\\
$[-1,-2,6,-7,-1]^*$ & $[1, 1, 1/2, 1/2, 1]$ & $[0, 1, 0, 5, 0]$  & $1$ & Yes\\
$[-7,5,-3,1,1]$ & $\rho/2$ & $\rho$  & $1$ & Yes\\
$[-5,3,-5,3,3]$ & $\rho/2$ & $\rho$  & $1$ & Yes\\
$[-7,5,-1,-1,-1]$ & $[1/2, 1/2, 1, 1, 1]$ & $[5, 1, 0, 0, 0]$  & $1$ & Yes\\
$[-2,4,-5,-2,4]^*$ & $[1, 1/2, 1/2, 1, 1/2]$ & $[0, 2, 0, 3, 0]$  & $1$ & Yes\\
$[-1,-5,3,-1,-1]$ & $[1, 1/2, 1/2, 1, 1]$ & $[3, 2, 0, 0, 0]$  & $1$ & Yes\\
$[-1,-1,-3,1,1]$ & $[1, 1, 1/2, 1/2, 1/2]$ & $[1, 3, 0, 0, 0]$  & $1$ & Yes\\
$[1,-3,1,-2,-2]$ & $[1/2, 1/2, 1/2, 1/2, 1]^*$ & $[0, 2, 0, 1, 2]$  & $1$ & Yes\\
$-\rho$ & $\rho/2$ & $\rho$  & $1$ & Yes\\
$-\rho$ & $\rho$ & $[0,0,0,0,0]$  & $1$ & Yes
\end{tabular}
\label{table-D5-scattered-part}
\end{table}

\section{The string part of $\widehat{E}_6^{\mathrm{d}}$}\label{sec-E6-string}

From now on, we set $G$ to be complex $E_6$, whose Dynkin diagram is in Figure 1, see page 687 of Knapp \cite{Kn} for more details. In particular, we note that $-w_0$ interchages $\alpha_1$ and $\alpha_6$, $\alpha_3$ and $\alpha_5$, while preserving $\alpha_2$, $\alpha_4$.

In this section, we use $\lambda=[\lambda_1, \dots, \lambda_6]$ to denote the weight $\sum_i\lambda_i\varpi_i$, where each $\lambda_i$ runs over $\frac{1}{2}\mathbb{P}$.  By Theorem \ref{thm-finite}, the string part of $\widehat{E}_6^{\mathrm{d}}$ comes from the scattered parts of $\widehat{L}_{\mathrm{ss}}^{\mathrm{d}}$, where $L\supseteq HA$ runs over all the proper $\theta$-stable Levi subgroups of $G$ and $L_{\mathrm{ss}}$ is its semisimple factor. Therefore, we can obtain the string part of $\widehat{E}_6^{\mathrm{d}}$ from the information in the previous section. Let us illustrate the process via examples.

\begin{example}\label{exam-HA}
Consider one extreme case $L=HA$. Then the corresponding string in $\widehat{E}_6^{\mathrm{d}}$ is:
$$
[\lambda_1, \dots, \lambda_6]\in \Lambda(e),
$$
 where each $\lambda_i\in\frac{1}{2}\mathbb{P}$. All of them are tempered representations. \hfill\qed
\end{example}

\begin{example}\label{exam-A2}
Consider the model representation of $A_2$, where $s=s_1 s_2 s_1$ and $\lambda=[1/2,1/2]$. This representation has $T$-parameter $[0,0]$ and $A$-parameter $[1,1]$. Recall that $s\rho(A_2)=[-1,-1]$.

There are five $\theta$-stable Levi subgroups of $G$ whose semisimple factors are of type $A_2$. The corresponding five strings in  $\widehat{E}_6^{\mathrm{d}}$ are listed as follows:
\begin{align*}
&[1/2, \lambda_2, 1/2, \lambda_4, \lambda_5, \lambda_6]\in \Lambda(s_1 s_3 s_1); \quad  [\lambda_1, \lambda_2, 1/2, 1/2, \lambda_5, \lambda_6]\in\Lambda (s_3 s_4 s_3)\\
&[\lambda_1, 1/2, \lambda_3, 1/2, \lambda_5, \lambda_6]\in\Lambda(s_2 s_4 s_2);\quad  [\lambda_1, \lambda_2, \lambda_3, 1/2, 1/2, \lambda_6]\in\Lambda(s_4 s_5 s_4)\\
&[\lambda_1, \lambda_2, \lambda_3, \lambda_4, 1/2, 1/2]\in\Lambda(s_5 s_6 s_5),
\end{align*}
 where each $\lambda_i\in\frac{1}{2}\mathbb{P}$.
We list the $T$-parameters and $A$-parameters for the first string below:
$$
[0, 2 \lambda_2, 0, 2 \lambda_4+1, 2 \lambda_5, 2\lambda_6], \quad [1, 0, 1, -1, 0, 0].
$$
\hfill\qed
\end{example}

\begin{example}\label{exam-D5}
Let us consider another extreme case. Namely, now $L_{\mathrm{ss}}$ is $D_5$. We focus on the first representation in Table \ref{table-D5-scattered-part}, where $s=s_1 s_3 s_2 s_1 s_4 s_5 s_3$ and $\lambda=[1,1/2,1/2,1/2,1/2]$.
This representation has $T$-parameter $[0, 3, -3, 3, 3]$ and $A$-parameter $[2, -2, 4, -2, -2]$. Recall that $s\rho(D_5)=[-2,5,-6, 4, 4]$.

There are two $\theta$-stable Levi subgroups of $G$ whose semisimple factors are of type $D_5$. The first one corresponds to $\{\alpha_1, \dots, \alpha_5\}$, while the second one corresponds to $\{\alpha_2, \dots, \alpha_6\}$. Let us denote the counterparts of the involution $s$ in them by $s'$ and $s''$, respectively. One can find that
$$
s'=s_4s_5s_1s_3s_2s_4s_1, \quad s''=s_4 s_5 s_6 s_5 s_3 s_2 s_4,
$$
and that
$$
s'\rho(E_6)=[-2,4,5,-6,4,\textbf{3}], \quad s''\rho(E_6)=[\textbf{3},4,4,-6,5,2].
$$
In each case, all the non-bolded coordinates  come from those of $s\rho(D_5)$ via permutations. Now the representation $J(\lambda, -s\lambda)$ of $D_5$ gives the following two strings of $\widehat{E}_6^{\mathrm{d}}$:
$$
[1,1/2,1/2,1/2,1/2, \lambda_6]\in\Lambda(s');\quad [\lambda_1, 1/2,1/2,1/2,1/2,1]\in \Lambda( s''),
$$
where $\lambda_1, \lambda_6\in\frac{1}{2}\mathbb{P}$.
We list their $T$-parameters and $A$-parameters below. They are
$$
[0, 3, 3, -3, 3, 2 \lambda_6+1], \quad [2, -2, -2, 4, -2, -1]
$$
for the first string, and
$$
[2\lambda_1+1, 3, 3, -3, 3, 0], \quad [-1, -2, -2, 4, -2, 2]
$$
for the second string. \hfill\qed
\end{example}

One can obtain other strings of $\widehat{E}_6^{\mathrm{d}}$ from the tables in Section \ref{sec-classical} without much difficulty. We omit this part. Instead, let us count the total number of strings in $\widehat{E}_6^{\mathrm{d}}$. Let $N_i$ be the number of $I$-strings in $\widehat{E}_6^{\mathrm{d}}$ such that $|I|=6-i$. Then $N_0 =1$  (see Example \ref{exam-HA}), and
\begin{align*}
N_1 &=6 N_{A_1}=6\\
N_2 &=5 N_{A_2}+10 N_{A_1} N_{A_1}=20,\\
N_3 &=5 N_{A_3} +10 N_{A_1} N_{A_2} + 5 N_{A_1} N_{A_1} N_{A_1}=45,\\
N_4 &= N_{D_4} + 4 N_{A_4}+4 N_{A_1}N_{A_3}+5 N_{A_1} N_{A_1} N_{A_2}+  N_{A_2} N_{A_2}=71,\\
N_5 &=2 N_{D_5}+N_{A_5}+2 N_{A_1} N_{A_4}+ N_{A_1} N_{A_2} N_{A_2}=70.
\end{align*}
Here recall that $N_G$ is the cardinality of the scattered part of $\widehat{G}^{\mathrm{d}}$. Therefore, $\widehat{E}_6^{\mathrm{d}}$ contains $\sum_{i=0}^{5} N_i=213$ strings in total.

\section{The scattered part of $\widehat{E}_6^{\mathrm{d}}$}

This section aims to report the scattered part of $\widehat{E}_6^{\mathrm{d}}$ using the computing method in Section \ref{sec-comp}. According to Theorem \ref{thm-finite}, we should focus on these $s$-families such that $I(s)$ is empty. That is, we should study those $s$-families such that any reduced expression of $s$ contains $s_1, \dots, s_6$.
There are $571$ such involutions, among which sit $103$ self-dual ones. Thus by using the $-w_0$ automorphism in Section \ref{sec-w0}, it boils down to consider $337$ $s$-families.

Let us provide a few more examples.

\begin{example}\label{exam-E6-2}
Consider the  involution $$s=s_6s_2s_4s_5s_3s_4s_1s_3s_2s_4s_5s_6s_3s_4s_1s_3s_2s_4s_1s_3s_2s_1.$$
 Note that $s\rho=[-1, -2, -1, -1, 10, -11]$ and $s$ is dual to the involution
 $$
s^\prime= s_6s_5s_4s_1s_3s_2s_4s_5s_6s_5s_4s_3s_2s_4s_5s_4s_3s_2s_4s_3s_2s_1.
 $$
 Note that $s^\prime\rho=[-11, -2, 10, -1, -1, -1]$ (recall Lemma \ref{lemma-involution-dual}).

 The first step of Section \ref{sec-comp} leaves us with $2475$ candidate representations. After carrying out the second step of Section \ref{sec-comp}, $35$ of them survive. By using \texttt{atlas}, we find that only $\lambda=[1,1,1,1, 1/2,1/2]$ gives a unitary representation $J(\lambda, -s\lambda)$, which has $T$-parameter
$[0, 0, 0, 0, 9, -9]$ and $A$-parameter $[2, 2, 2, 2, -8, 10]$. By looking at its $K$-types, we know that it has a unique spin-lowest $K$-type $[0, 1, 0, 0, 0, 9]$ which occurs once.  Moreover,
 $$
 \|[0, 1, 0, 0, 0, 9]\|_{\mathrm{spin}}=\|2\lambda\|.
 $$
 Thus it is a member of $\widehat{E}_6^{\mathrm{d}}$ by Proposition \ref{prop-D-spin-lowest}.
 This gives the $13$th row of Table \ref{table-E6-scattered-part}.\hfill\qed
\end{example}

\begin{table}
\centering
\caption{The scattered part of $\widehat{E}_6^{\mathrm{d}}$ (folded version)}
\begin{tabular}{l|c|c|c|r}
$s\rho$ &   $\lambda$   & spin LKT & mult &  u-small\\
\hline
$[-2,5,6,-7,6,-2]$ & $[1, 1/2, 1/2, 1/2, 1/2, 1]$ & $[1, 1, 0, 3, 0, 1]$  & $1$ & Yes\\
$[-4,-2,3,6,-8,6]^*$ & $[1/2, 1, 1/2, 1/2, 1/2, 1/2]$ & $[2, 1, 0, 1, 2, 1]$  & $1$ & Yes\\
$[-5,-7,3,5,3,-5]$ & $\rho/2$ & $\rho$  & $1$ & Yes\\
$[-1,-1,-2,8,-9,7]^*$ & $[1, 1, 1, 1/2, 1/2, 1/2]$ & $[0, 0, 1, 0, 5, 1]$  & $1$ & Yes\\
$[-2,-8,-2,7,4,-5]^*$ & $[1, 1/2, 1, 1/2, 1/2, 1/2]$ & $[1, 4, 1, 0, 0, 3]$  & $1$ & Yes\\
$[-8,-1,6,-1,6,-8]$ & $[1/2, 1, 1/2, 1, 1/2, 1/2]$ & $[4, 0, 1, 0, 1, 4]$  & $1$ & Yes\\
$[-10,5,8,-6,4,-2]^*$ & $[1/2, 1/2, 1/2, 1/2, 1/2, 1]$ & $[1, 0, 2, 1, 1, 1]$  & $1$ & Yes\\
$[-1,-11,-1,9,-1,-1]$ & $[1, 1/2, 1, 1/2, 1, 1]$ & $[0, 7, 0, 1, 0, 0]$  & $1$ & Yes\\
$[-8,7,6,-8,6,-2]^*$ & $[1/2, 1/2, 1/2, 1/2, 1/2, 1]$ & $[1, 0, 2, 1, 1, 1]$  & $1$ & Yes\\
$[-1,-1,-6,4,7,-9]^*$ & $[1, 1, 1/2, 1/2, 1/2, 1/2]$ & $[2, 0, 2, 0, 1, 3]$  & $1$ & Yes\\
$[-1,-1,8,-10,8,-1]$ & $[1, 1, 1/2, 1/2, 1/2, 1]$ & $[0, 4, 0, 2, 0, 0]$  &  $1$& Yes \\
$[-11,-3,9,1,-3,1]^*$ & $\rho/2$ & $\rho$  & $1$ & Yes\\
$[-1,-2,-1,-1,10,-11]^*$ & $[1, 1, 1, 1, 1/2, 1/2]$ &  $[0, 1, 0, 0, 0, 9]$ & $1$ & Yes \\
$[6,-1,-8,6,-8,6]$ & $[1/2, 1, 1/2, 1/2, 1/2, 1/2]$ &  $[1, 2, 0, 2, 0, 1]$ & $1$ & Yes \\
$[-1,-1,-10,8,-1,-1]^*$ & $[1, 1, 1/2, 1/2, 1, 1]$ & $[7, 2, 0, 0, 0, 0]$   & $1$ & Yes\\
$[-2,7,-1,-8,6,-1]^*$ & $[1, 1/2, 1, 1/2, 1/2, 1]$ &  $[5, 3, 0, 0, 0, 0]$ & $1$ & Yes \\
$[-5,3,3,-5,3,-5]$ & $\rho/2$ & $\rho$  & $1$ & Yes \\
$[-2,-6,-1,4,-5,3]^*$ & $[1, 1/2, 1, 1/2, 1/2, 1/2]$ &  $[4, 3, 0, 0, 0, 1]$ & $1$ & Yes \\
$[-2,-1,1,-3,1,-2]$ & $[1/2, 1, 1/2, 1/2, 1/2, 1]^*$ & $[2, 3, 0, 0, 0, 3]$   & $1$  & Yes \\
$-\rho$ & $\rho/2$ & $\rho$  & $1$ & Yes \\
$-\rho$ & $\rho$ & $[0,0,0,0,0,0]$  & $1$ & Yes
\end{tabular}
\label{table-E6-scattered-part}
\end{table}

\begin{example}\label{exam-E6-3}
Consider the  involution $$s=s_6s_5s_3s_4s_1s_3s_2s_4s_5s_6s_5s_4s_3s_2s_4s_5s_4s_3s_2s_4s_1s_3s_2s_1.$$
 Note that $s\rho=[-1, -1, -10, 8, -1, -1]$ and $s$ is dual to the involution
 $$
s^\prime= s_5s_6s_2s_4s_5s_3s_4s_1s_3s_2s_4s_5s_6s_3s_4s_1s_3s_2s_4s_5s_1s_3s_2s_1.
 $$
 Note that $s^\prime\rho=[-1, -1, -1, 8, -10, -1]$ (recall Lemma \ref{lemma-involution-dual}).

 The first step of Section \ref{sec-comp} leaves us with $1145$ candidate representations. After carrying out the second step of Section \ref{sec-comp}, $17$ of them survive. By using \texttt{atlas}, we find that only $\lambda=[1,1, 1/2,1/2, 1,1]$ gives a unitary representation $J(\lambda, -s\lambda)$, which has $T$-parameter
$[0, 0, -7, 7, 0, 0]$ and $A$-parameter $[2, 2, 8, -6, 2, 2]$. By looking at its $K$-types, we know that it has a unique spin-lowest $K$-type $[7,2,0,0,0,0]$ which occurs once. Moreover,
 $$
 \|[7,2,0,0,0,0]\|_{\mathrm{spin}}=\|2\lambda\|.
 $$
 Thus it is a member of $\widehat{E}_6^{\mathrm{d}}$ by Proposition \ref{prop-D-spin-lowest}.
 This gives the $15$th row of Table \ref{table-E6-scattered-part}.\hfill\qed
\end{example}

The final result is given in Table \ref{table-E6-scattered-part}. Note that the representation in the last row is the trivial one. The scattered part of $\widehat{E}_6^{\mathrm{d}}$ consists of $33$ representations in total. That is, $N_{E_6}=33$.

\end{document}